\documentclass{ifacconf}

\usepackage{graphicx}      
\usepackage{natbib}        

\usepackage{amsmath} 
\usepackage{amssymb}  
\usepackage{epstopdf}
\usepackage{subfigure,color,balance}
\usepackage{verbatim}
\usepackage{cases}
\usepackage{enumerate}
\usepackage{balance}
\usepackage{mathbbol}
\usepackage{framed}

\newtheorem{remark}{Remark}
\newtheorem{problem}{Problem}
\newtheorem{definition}{Definition}

\begin{document}
\begin{frontmatter}

\title{Optimal Formation of Autonomous Vehicles \\ in Mixed Traffic Flow\thanksref{footnoteinfo}}

\thanks[footnoteinfo]{This work is supported by National Key R\&D Program of China with 2016YFB0100906. K. Li and J. Wang contributed equally to this work. All correspondence should be sent to Y. Zheng. }

\author[First]{Keqiang Li}
\author[First]{Jiawei Wang}
\author[Second]{Yang Zheng}

\address[First]{School of Vehicle and Mobility, Tsinghua University,
	China \\(e-mail: likq@tsinghua.edu.cn, wang-jw18@mails.tsinghua.edu.cn)}
\address[Second]{School of Engineering and Applied Sciences, Harvard University,
	USA (e-mail: zhengy@g.harvard.edu)}

\begin{abstract}                
Platooning of multiple autonomous vehicles has attracted significant attention in both academia and industry. Despite its great potential, platooning is not the only choice for the formation of autonomous vehicles in mixed traffic flow, where autonomous vehicles and human-driven vehicles (HDVs) coexist. In this paper, we investigate the optimal formation of autonomous vehicles that can achieve an optimal system-wide performance in mixed traffic flow. Specifically, we consider the optimal $\mathcal{H}_2$ performance of the entire traffic flow, reflecting the potential of autonomous vehicles in mitigating traffic perturbations. Then, we formulate the optimal formation problem as a set function optimization problem. Numerical results reveal two predominant optimal formations: uniform distribution and platoon formation, depending on traffic parameters. In addition, we show that 1) the prevailing platoon formation is not always the optimal choice; 2) platoon formation might be the worst choice when HDVs have a poor string stability behavior. These results suggest more opportunities for the formation of autonomous vehicles, beyond platooning, in mixed traffic flow.
\end{abstract}


\end{frontmatter}

\section{Introduction}
For a series of human-driven vehicles (HDVs), it is known that small perturbations may be accumulated and amplified, finally leading to stop-and-go waves. This phenomenon, also known as phantom traffic jam, has resulted in a great loss of travel efficiency and fuel economy \citep{sugiyama2008traffic}. The emergence of autonomous vehicles is expected to smooth traffic flow and improve traffic efficiency significantly. In particular, instead of controlling each vehicle separately, cooperative formation and control of multiple autonomous vehicles will revolutionize future transportation systems \citep{li2017dynamical}.

One typical example of cooperative formation is vehicle platooning, which has attracted increasing attention in the past decades. In a platoon, adjacent autonomous vehicles are regulated to maintain the same desired velocity while keeping a pre-specified inter-vehicle distance. Both rigorous theoretical analysis (e.g., \citealp{zheng2015stability}) and real-world experiments (e.g., \citealp{milanes2013cooperative}) have confirmed the great potential of vehicle platooning in achieving higher traffic efficiency, better driving safety, and lower fuel consumption. However, platooning technologies typically require all the involved vehicles to have autonomous capabilities. As the gradual deployment of autonomous vehicles in practice, there will be a long transition phase of mixed traffic flow, where HDVs and autonomous vehicles coexist. This brings a challenge for the practical implementation of vehicle platooning. Since autonomous vehicles are usually distributed randomly in real traffic flow, several maneuvers, including joining, leaving, merging, and splitting, need to be performed to form the neighboring autonomous vehicles into a platoon \citep{amoozadeh2015platoon}. Despite the great potential of platooning, recent works have revealed the possible negative impacts of these maneuvers, which might even cause undesired congestions; see, e.g., \cite{mena2018impact}. These results suggest that forming a platoon for autonomous vehicles might not be necessary in mixed traffic flow.

\begin{figure}[t]
	\centering
	\subfigure[Uniform distribution]
	{
		\label{Fig:Formationillustration_Unifrom}
		\includegraphics[scale=0.2]{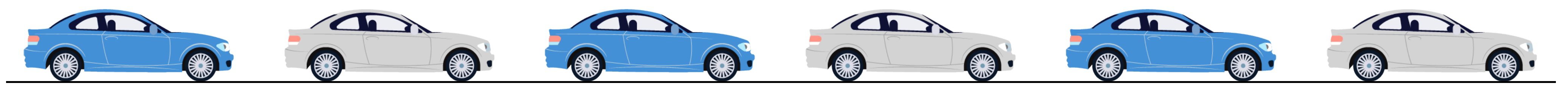}
	}
	\subfigure[Random formation]
{ \label{Fig:Formationillustration_Random}
	\includegraphics[scale=0.2]{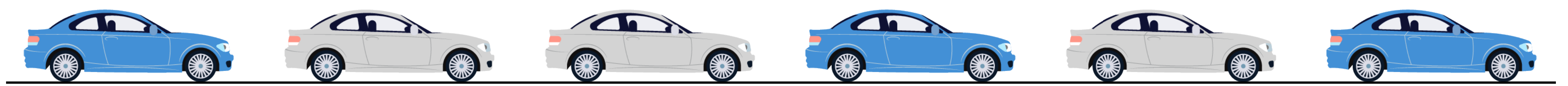}
}
	\subfigure[Platoon formation]
	{\label{Fig:Formationillustration_Platoon}
		\includegraphics[scale=0.2]{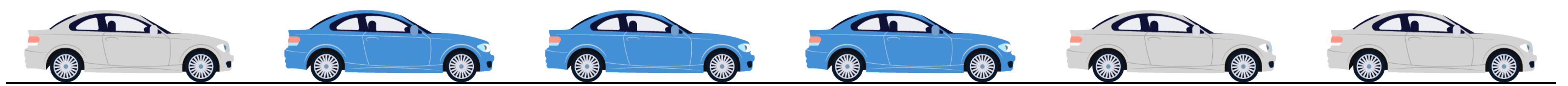}
	}
	\vspace{-3mm}
	\caption{Possible formations of autonomous vehicles in mixed traffic flow. Gray vehicles and blue vehicles denote HDVs and autonomous vehicles, respectively.}
	\label{Fig:FormationIllustration}
\end{figure}

In fact, platooning is not the only formation option for autonomous vehicles in mixed traffic flow. Possible choices can be much more diverse, since any combination form of HDVs and autonomous vehicles is feasible. For example, uniform distribution of autonomous vehicles (see Fig.~\ref{Fig:Formationillustration_Unifrom} for illustration) could be another simple formation. It is of great importance to investigate the formation for autonomous vehicles that can achieve optimal system-wide performance for the entire mixed traffic flow. However, most existing works on mixed traffic flow focus on the influence of penetration rates; see, e.g., \cite{talebpour2016influence}. These works typically assume either a random formation (Fig.~\ref{Fig:Formationillustration_Random}) or the platoon formation (Fig.~\ref{Fig:Formationillustration_Platoon}). The role of different formations of autonomous vehicles in mixed traffic flow has not been well-understood. In particular, whether the prevailing platoon formation performs better than other formations in mixed traffic flow remains unclear. In this paper, our main focus is to identify an optimal formation for autonomous vehicles in  mixed traffic environment.

To address this problem, we utilize the notion of \emph{Lagrangian control} of traffic flow \citep{stern2018dissipation}, where autonomous vehicles serve as \emph{mobile actuators} for traffic control. Its effectiveness in dampening traffic waves has been validated in the case of one single autonomous vehicle, including theoretical analysis \citep{cui2017stabilizing,zheng2020smoothing}, small-scale real-world experiments \citep{stern2018dissipation} and large-scale numerical simulations \citep{vinitsky2018lagrangian}. Along this direction, it is important to investigate the influence of different placements of multiple autonomous vehicles on mixed traffic systems. A closely related topic is the so-called \emph{actuator placement} problem, which has been discussed in many dynamical systems, including mechanical systems \citep{hiramoto2000optimal}, power grids \citep{qin2018submodularity} and biological networks \citep{gu2015controllability}. It has also attracted extensive attention to investigate an optimal placement of multiple actuators in order to maximize certain system performance~\citep{olshevsky2014minimal,summers2015submodularity}. To the best of our knowledge, the placement of autonomous vehicles in traffic flow has not been discussed in the literature, and the previous results on actuator placement~\citep{olshevsky2014minimal,summers2015submodularity} are not directly applicable.

In this paper, we focus on the problem of optimal formation of autonomous vehicles in a ring-road mixed traffic system. The formation of autonomous vehicles is characterized by their placement pattern, i.e., their locations, in traffic flow. Particularly, we aim to answer whether platooning achieves a better performance than other formations in mixed traffic flow. We formulate this problem as maximization of a formation value function, which is a set function representing the performance of the global traffic system. Specifically, the $\mathcal{H}_2$ performance is considered, and submodularity is discussed for this formulation. The contributions of this paper are as follows.

\begin{enumerate}
	\item
	A set function optimization formulation is proposed to model the optimal formation problem of autonomous vehicles in mixed traffic flow. A global $\mathcal{H}_2$ optimal controller is considered for the autonomous vehicles, which can reveal the maximum potential of a given formation of autonomous vehicles in mitigating traffic perturbations. We show that this optimization problem is in general not submodular.
	\item
	We show that platooning is not always the optimal formation for autonomous vehicles in mixed traffic flow. In fact, numerical results reveal two predominant optimal formations: platoon formation and uniform distribution. Furthermore, we find that the optimal formation relies heavily on the string stability performance of the HDV car-following behavior. When HDVs have a poor string stability behavior, platoon formation might be the worst choice.
\end{enumerate}

The rest of this paper is organized as follows. Section~\ref{Section:Modeling} introduces the modeling for mixed traffic systems and the problem statement. Formulation and analysis of the optimal formation problem is presented in Section~\ref{Section:OptimalFormulation}. Section~\ref{Section:Results} shows numerical results, and we conclude this paper in Section~\ref{Section:Conclusion}.

\section{SYSTEM MODELING AND PROBLEM STATEMENT}

\label{Section:Modeling}

In this section, we present a dynamic model of mixed traffic systems with multiple autonomous vehicles in a ring-road setup, and introduce the optimal formation problem. 

\begin{figure}[t]
	\centering
		\includegraphics[scale=0.325]{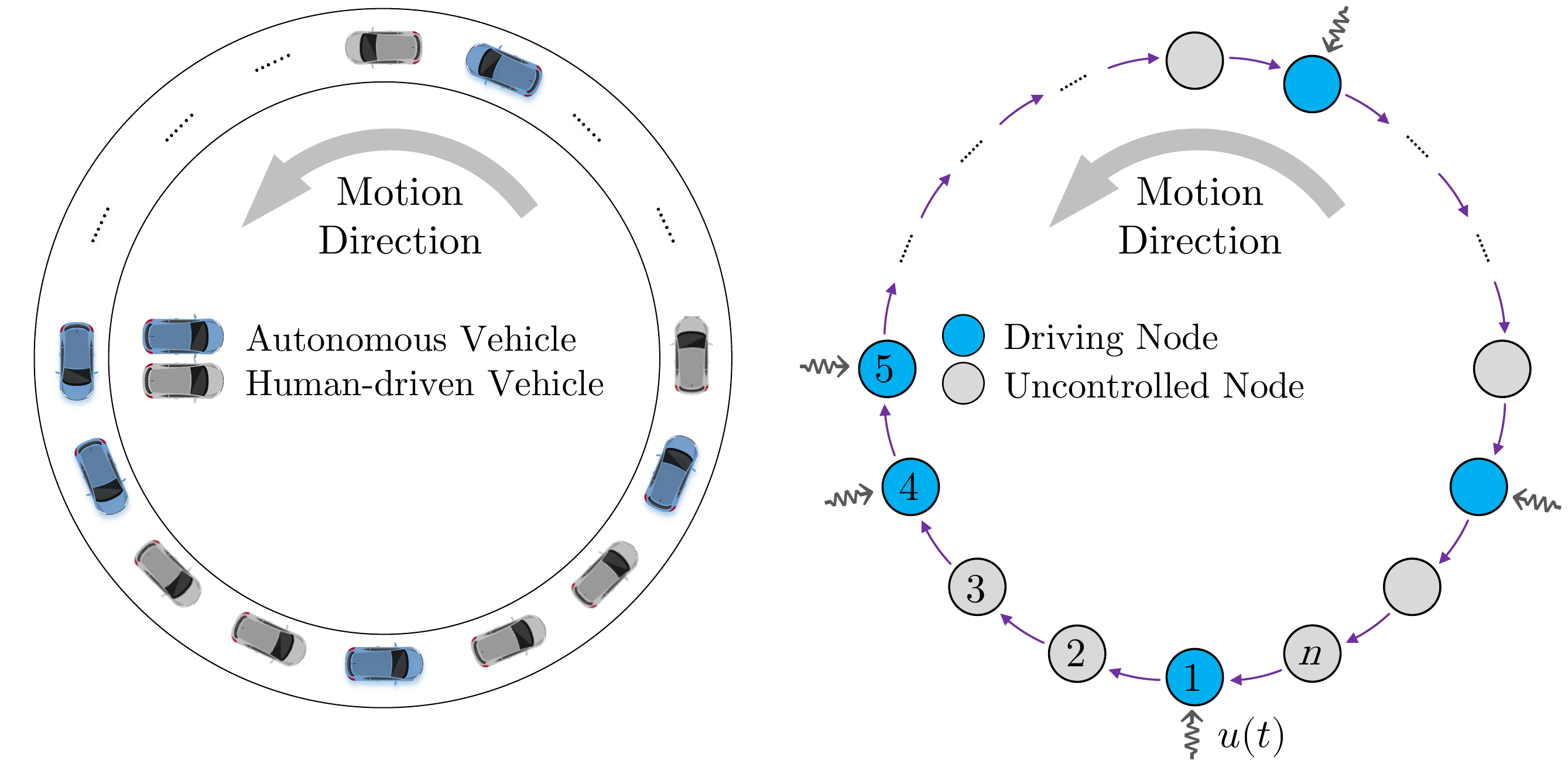}
	\vspace{-3mm}
	\caption{System model. (a) The single-lane ring road scenario with autonomous vehicles and human-driven vehicles. (b) A simplified network system schematic. }
	\label{Fig:SystemModel}
\end{figure}

\subsection{Modeling Mixed Traffic Systems}

As shown in Fig.~\ref{Fig:SystemModel}, we consider a single-lane ring road of length $L$ with $n$ vehicles, among which there are $k$ autonomous vehicles and $n-k$ human-driven vehicles. The ring road setup has been widely used in the literature~\citep{sugiyama2008traffic,cui2017stabilizing,stern2018dissipation,zheng2020smoothing}. This setup represents a simplified closed traffic system with no boundary conditions, and also corresponds to a straight road of infinite length and periodic traffic dynamics. 

The vehicles are indexed from $1$ to $n$, and we define $\Omega=\{1,2,\ldots,n\}$. We characterize the formation of autonomous vehicles by their placement pattern in mixed traffic flow, which is represented as a set $S=\{i_1,\ldots,i_k \}\subseteq\Omega$, where $i_1,\ldots,i_k$ denote the indices of autonomous vehicles. We denote the position, velocity and acceleration of vehicle $i$ as $p_i$, $v_i$ and $a_i$, respectively. The spacing of vehicle $i$, i.e., its relative distance from vehicle $i-1$, is defined as $s_i=p_{i-1}-p_i$. Then the relative velocity can be expressed as $\dot{s}_i=v_{i-1}-v_i$. The vehicle length is ignored without loss of generality.

Based on existing HDV models, the longitudinal dynamics of an HDV can be described by the following nonlinear process~\citep{orosz2010traffic,treiber2013traffic}
\begin{equation}\label{Eq:HDVNonlinearModel}
\dot{v}_i (t)=F(s_i (t),\dot{s}_i (t),v_i (t)),\;i\notin S,
\end{equation}
meaning that the acceleration of an HDV is determined by the relative distance, relative velocity and its own velocity. In equilibrium traffic flow, where $\dot{v}_i=0$ for $i=1,2,\ldots,n$, each vehicle moves with the same equilibrium velocity $v^*$ and spacing $s^*$. Around the equilibrium state $(s^*,v^*)$, we define the error state as
\begin{equation}
\tilde{s}_i(t)=s_i(t)-s^*, \;
\tilde{v}_i(t)=v_i(t)-v^*.
\end{equation}
Applying the first-order Taylor expansion to \eqref{Eq:HDVNonlinearModel}, we can derive a linearized model for each HDV around the equilibrium state
\begin{equation}\label{Eq:LinearHDVModel}
\begin{cases}
\dot{\tilde{s}}_i(t)=\tilde{v}_{i-1}(t)-\tilde{v}_i(t),\\
\dot{\tilde{v}}_i(t)=\alpha_{1}\tilde{s}_i(t)-\alpha_{2}\tilde{v}_i(t)+\alpha_{3}\tilde{v}_{i-1}(t),\\
\end{cases} i\notin S,
\end{equation}
with $\alpha_1 = \frac{\partial F}{\partial s}, \alpha_2 = \frac{\partial F}{\partial \dot{s}} - \frac{\partial F}{\partial v}, \alpha_3 = \frac{\partial F}{\partial \dot{s}}$  evaluated at the equilibrium state $(s^*,v^*)$. According to the real driving behavior, we have $\alpha _{1}>0$, $\alpha _{2}> \alpha _{3}>0$~\citep{jin2016optimal,cui2017stabilizing}.

For each autonomous vehicle, the acceleration signal is directly used as the control input $u_{i}(t),$ and its car-following model is thus given by 
\begin{equation}
\begin{cases}
\dot{\tilde{s}}_{i}(t)=\tilde{v}_{i-1}(t)-\tilde{v}_{i}(t),\\
\dot{\tilde{v}}_{i}(t)=u_{i}(t),
\end{cases} i \in S.
\end{equation}

To model traffic perturbations, we assume there exists a scalar disturbance signal $\omega_{i}(t)$ with finite energy in the acceleration of vehicle $i$ $(i \in\Omega) .$ Lumping all the error states into one global state $x(t)=\left[\tilde{s}_{1}(t), \ldots, \tilde{s}_{n}(t), \tilde{v}_{1}(t), \ldots, \tilde{v}_{n}(t)\right]^{T}$ and letting $\omega(t)=[\omega_{1}(t)$ $, \ldots, \omega_{n}(t)]^{T}$, $u(t)=\left[u_1(t),\ldots,u_n(t)\right]^T$, the state-space model for the mixed traffic system is then obtained
\begin{equation} \label{Eq:SystemModel}
	\dot{x}(t)=A_{S} x(t)+B_{S} u(t)+H \omega(t),
\end{equation}
where we have
\begin{equation*}
	\begin{aligned}
	A_{S}&=\begin{bmatrix}{0} & {M_{1}} \\ {\alpha_{1}\left(I_{n}-D_{S}\right)} & {P_{S}}\end{bmatrix} \in \mathbb{R}^{2 n \times 2 n}, \\ B_{S}&=\begin{bmatrix}{\mathbb{e}_{i_{1}}, \mathbb{e}_{i_{2}}, \ldots, \mathbb{e}_{i_{k}}}\end{bmatrix} \in \mathbb{R}^{2 n \times k}, \\
	H&=\begin{bmatrix}{0} \\ {I_{n}}\end{bmatrix}\in \mathbb{R}^{2 n \times n},
	\end{aligned}
\end{equation*}
and
\begin{equation*}
	\begin{aligned}
	M_1 &= \begin{bmatrix}
	-1 & & \cdots &1\\
	1&-1& &\\
	 & \ddots & \ddots &\\
	 & & 1 &-1 	\end{bmatrix}\in \mathbb{R}^{n\times n},\\
	D_S &= \mathrm{diag}\left(\delta_1,\delta_2,\ldots,\delta_n\right)\in \mathbb{R}^{n\times n},\\
	P_S &= \begin{bmatrix}
	-\alpha_2 \bar{\delta}_1 & & \cdots &\alpha_3 \bar{\delta}_1\\
	\alpha_3 \bar{\delta}_2&-\alpha_2 \bar{\delta}_2& &\\
	& \ddots & \ddots &\\
	& & \alpha_3 \bar{\delta}_n &-\alpha_2 \bar{\delta}_n 	\end{bmatrix}\in \mathbb{R}^{n\times n}.
	\end{aligned}
\end{equation*}
Throughout this paper, we use $I_{n}$ and $\mathrm{diag}(\cdot)$ to denote an identity matrix of size $n$ and a diagonal matrix, respectively. We also define a bool variable $\delta_{i}$ to indicate whether vehicle $i$ is an autonomous vehicle, i.e., $ \delta_{i}=0,$ if $i \notin S$; $\delta_{i}=1,$ if $i \in S$. Let $\bar{\delta}_{i}=1-\delta_{i}$ indicate whether vehicle $i$ is an HDV. In the input matrix $B_{S}$, the vector $\mathbb{e}_{i_r}$ is a $2 n \times 1$ unit vector $(r=1,2, \ldots, k)$, with the $\left(i_{r}+n\right)$-th entry being one and the others being zeros.

\begin{remark}
	Note that another mathematical model was introduced in \cite{zheng2020smoothing} to describe the dynamics of a ring-road mixed traffic system with multiple autonomous vehicles. Since the state vector therein can be transformed to $x(t)$ in \eqref{Eq:SystemModel} by a permutation matrix, model \eqref{Eq:SystemModel} is essentially equivalent to that in \cite{zheng2020smoothing}. In this paper, we choose the form \eqref{Eq:SystemModel} 
	due to its convenience 
	to reflect the relationship between system matrices and the formation $S$ of autonomous vehicles.
\end{remark}

\begin{remark}
It is demonstrated that autonomous vehicles can be utilized as mobile actuators for traffic control, leading to the notion of Lagrangian control of traffic flow~\citep{stern2018dissipation}. The potential of one single autonomous vehicle in stabilizing traffic flow and improving traffic velocity has been revealed in recent works~\citep{cui2017stabilizing,stern2018dissipation,zheng2020smoothing,wang2019controllability}. When multiple autonomous vehicles coexist, their specific formation plays a key role on the global performance of the entire traffic flow. The prevailing vehicle platooning is a straightforward choice for the formation (see Fig~\ref{Fig:FormationExample}(a)). However, whether platooning is the optimal one remains unclear.
\end{remark}

\subsection{Problem Statement}

In this paper, our main focus is to identify an optimal formation that maximizes a system-wide performance metric of the entire traffic system. 

\begin{problem} \label{Pr:Main}
	Assume there are $k$ autonomous vehicles in the ring-road mixed traffic system \eqref{Eq:SystemModel}. Find an optimal formation, i.e., $S=\left\{i_{1}, \ldots, i_{k}\right\} \subseteq \Omega,$ for the autonomous vehicles, which achieves the optimal system-wide performance for the entire traffic flow.
\end{problem}

\begin{figure}[t]
	\centering
	\subfigure[]
	{ \label{Fig:ExamplePlatoon}
		\includegraphics[scale=0.55]{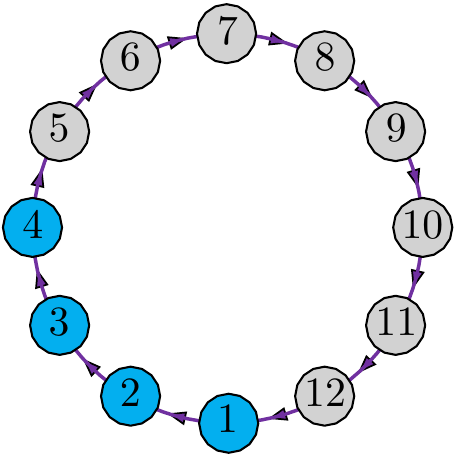}
	}
	\subfigure[]
	{\label{Fig:ExampleUniform}
		\includegraphics[scale=0.55]{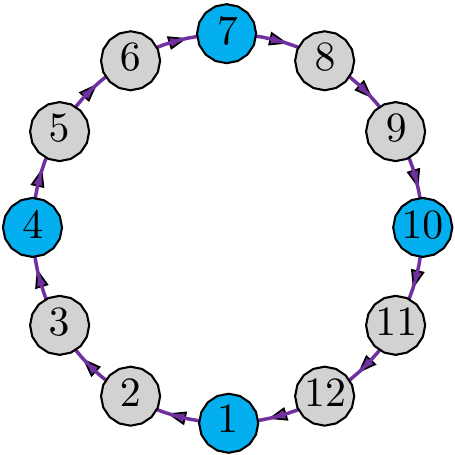}
	}
\subfigure[]
{
	\label{Fig:ExampleAbnormal}
	\includegraphics[scale=0.55]{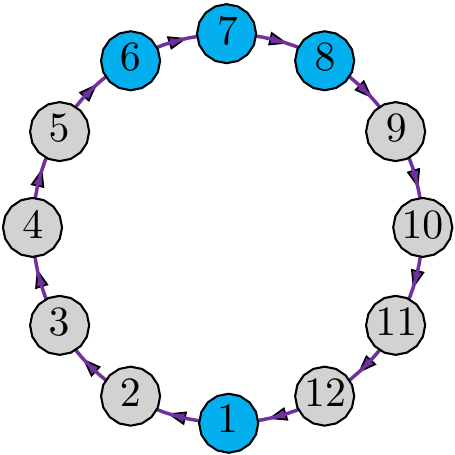}
}
	\vspace{-2mm}
	\caption{Possible formations when $n=12,k=4$. Gray nodes: HDVs; blue nodes: autonomous vehicles. (a) Platoon formation ($S=\{1,2,3,4\}$). (b) Uniform distribution ($S=\{1,4,7,10\}$). (c) Abnormal formation ($S=\{1,6,7,8\}$).}
	\label{Fig:FormationExample}
\end{figure}

In Fig.~\ref{Fig:FormationExample}, we illustrate three examples of the formation of autonomous vehicles in the ring-road mixed traffic system, when $n=12$, $k=4$. Possible formations include platoon formation (Fig.~\ref{Fig:ExamplePlatoon}), uniform distribution (Fig.~\ref{Fig:ExampleUniform}) and other abnormal cases (Fig.~\ref{Fig:ExampleAbnormal}). We are interested in whether the prevailing platoon formation is the optimal choice for the mixed traffic scenario.
To formulate Problem \ref{Pr:Main} mathematically, we utilize the formation $S$ as the decision variable. To quantify the performance of certain formation, we introduce a formation value function, $J(S): 2^{\Omega} \rightarrow \mathbb{R},$ to measure the system-wide 
performance for a given formation $S \subseteq \Omega$. Note that $J(S)$ is a set function, and a larger value of $J(S)$ indicates a better performance. The cardinality of $S$ is denoted as $|S|$. Then Problem 1 can be formulated abstractly as follows
\begin{equation} \label{Eq:ProblemFormulation}
\begin{aligned}
	\max_S  \quad &J(S)\\
\mathrm{s.t.}\quad  & S \subseteq \Omega,|S|=k
\end{aligned}
\end{equation}
where the optimal solution $S^*$ offers the optimal formation for autonomous vehicles in mixed traffic flow.

\begin{remark}
	Formulation~\eqref{Eq:ProblemFormulation} is a standard set function optimization problem, which has been widely used in the problem of actuator placement in dynamical systems; see, e.g., \cite{qin2018submodularity}. For a linear time-invariant system given by $\dot{x}=Ax+Bu$, most existing results typically consider the case where the placement decision only determines the input matrix $B$~\citep{olshevsky2014minimal,summers2015submodularity}. In mixed traffic flow, however, autonomous vehicles and HDVs have distinct dynamics. When we choose a different formation for autonomous vehicles, the system matrix $A$ will also be changed. Therefore, in our system model \eqref{Eq:SystemModel}, both the system matrix $A_S$ and input matrix $B_S$ rely on the formation $S$ of autonomous vehicles. The results in \cite{olshevsky2014minimal} and \cite{summers2015submodularity} are not directly applicable.
\end{remark}

\begin{remark}
	Note that the choice of the system performance metric $J(S)$ has a great influence on the optimal solution to~\eqref{Eq:ProblemFormulation}. In previous works on optimal actuator placement, controllability criteria have received the most attention~\citep{olshevsky2014minimal,summers2015submodularity}. However, it has been shown in \cite{zheng2020smoothing} that a ring-road mixed traffic system is not completely controllable but is stabilizable, when $\vert S \vert \geq 1$. In this paper, we consider an  $\mathcal{H}_2$ optimal control performance metric, which can reveal the maximum potential of a given formation of autonomous vehicles in mitigating traffic perturbations. The detailed formulation is presented in the next section.
\end{remark}



\section{Problem Formulation and Analysis}

\label{Section:OptimalFormulation}

In this section, a global $\mathcal{H}_2$ optimal controller is applied to the autonomous vehicles and the resulting optimal $\mathcal{H}_2$ performance value is chosen as the explicit form of the formation value function $J(S)$. Submodularity of this specific formulation is also discussed.

\subsection{Global Optimal Controller}
For a given formation $S$ of autonomous vehicles, we consider a global state feedback controller
\begin{equation}
	u=-K_{S} x, \;K_{S} \in \mathbb{R}^{2 n \times k}.
\end{equation}
The control objective is to achieve an optimal performance for the global mixed traffic system via controlling the autonomous vehicles. Specifically, the control target is to minimize the influence of undesired perturbations $\omega(t)$ on the mixed traffic system. Note that the optimal feedback gain $K_{S}$ relies on the choice of the formation $S$.

We use $z(t)$ to denote a performance output for the global mixed traffic system
 \begin{equation} \label{Eq:Output}
z(t) = \begin{bmatrix} Q^{\frac{1}{2}} \\0 \end{bmatrix}x(t) +  \begin{bmatrix} 0 \\R^{\frac{1}{2}} \end{bmatrix}u(t),
\end{equation}
where $Q^{\frac{1}{2}} = \mathrm{diag}\left(\gamma_s,\ldots,\gamma_s,\gamma_v,\ldots,\gamma_v\right) $ and $R^{\frac{1}{2}} = \mathrm{diag}(\gamma_u,$
$\ldots,\gamma_u) $. The weight coefficients $\gamma_s,\gamma_v,\gamma_u>0$ represent the penalty for spacing error, velocity error and control input, respectively. When applying the controller $u=-K_S x$, the dynamics of the closed-loop mixed traffic system then become
\begin{equation}
\begin{aligned}
\dot{x}(t) &= (A_S - B_S K_S)x(t) + Hw(t), \\
z(t) &= \begin{bmatrix} Q^{\frac{1}{2}} \\ -R^{\frac{1}{2}}K_S \end{bmatrix}x(t).
\end{aligned}
\end{equation}
We use the $\mathcal{H}_2$ norm of the transfer function $G_S$ from disturbance $\omega$ to output $z$ to describe the influence of perturbations on the traffic system. Then the optimal feedback gain $K_S$ of the autonomous vehicles 
can be obtained by solving the following optimization problem
\begin{equation} \label{Eq:H2Control}
	\min_{K_S} \; \lVert G_S \rVert_2^2
\end{equation}
where $\lVert \cdot \rVert_2$ denotes the $\mathcal{H}_2$ norm. The optimization problem \eqref{Eq:H2Control} is in the standard form of $\mathcal{H}_2$ optimal control~\citep{skogestad2007multivariable}. Here, we directly present a convex reformulation for \eqref{Eq:H2Control} as follows~\citep{zheng2020smoothing}
  \begin{equation}\label{Eq:LMIOptimalControl}
\begin{aligned}
\min_{X,Y,Z} \; & {\mathrm{Tr}}(QX)+{\mathrm{Tr}}\left(RY\right) \\
\mathrm{s.t.} \; & (A_S X-B_S Z)+(A_S X-B_S Z)^{T} + HH^{T} \preceq 0, \\
&\begin{bmatrix}
Y & Z \\ Z^T & X\end{bmatrix}\succeq 0,\; X \succ 0,
\end{aligned}
\end{equation}
where $\mathrm{Tr}(\cdot)$ denotes the trace of a matrix. Problem \eqref{Eq:LMIOptimalControl} can be reformulated into a standard semidefinite program, which can be solved efficiently via existing solvers, e.g., Mosek~\citep{mosek2010mosek}.

\subsection{Reformulation of Optimal Formation}
For a given formation decision $S$, the optimal feedback gain $K_S$ can be obtained by solving \eqref{Eq:LMIOptimalControl}. Meanwhile, the optimal value of $\min_{K_S}  \lVert G_S \rVert_2^2 $ indicates the minimum influence of perturbations on the entire traffic flow when the autonomous vehicles are optimally controlled. Accordingly, the specific expression of the formation value function $J(S)$ in \eqref{Eq:ProblemFormulation} can be given by
\begin{equation} \label{Eq:H2ObjectiveValue}
	J(S):=-\min_{K_S}  \lVert G_S \rVert_2^2.
\end{equation}
The negative sign exists for normalization. 
Now, the original problem \eqref{Eq:ProblemFormulation} of optimal formation of autonomous vehicles in mixed traffic flow can be reformulated as
\begin{equation} \label{Eq:H2ProblemFormulation}
\begin{aligned}
\max_{S}  \quad &J(S)=-\min_{K_S}  \lVert G_S \rVert_2^2\\
\mathrm{s.t.}\quad  &S \subseteq \Omega,|S|=k
\end{aligned}
\end{equation}
In \eqref{Eq:H2ProblemFormulation}, the optimization problem \eqref{Eq:LMIOptimalControl} needs to be solved first to calculate the specific value of $J(S)$ for a given formation decision $S$. Since it is proved in \cite{zheng2020smoothing} that the mixed traffic system with one or more autonomous vehicles is always stabilizable, there exist stabilizing feedback gains $K_S$ under which the $\mathcal{H}_2$ norm of $G_S$ is finite, when $|S|\geq 1$.

\subsection{Submodularity Analysis}
Based on combinatorial optimization problem~\eqref{Eq:H2ProblemFormulation}, we can obtain the optimal formation solution by enumerating all possible subsets of cardinality $k$. 
Here, we attempt to investigate whether \eqref{Eq:H2ProblemFormulation} possesses certain useful properties that lead to practically efficient algorithms. In particular, we consider the property of submodularity, which plays a significant role in set function optimization problems. For submodular functions, a simple greedy algorithm can return a near-optimal solution~\citep{nemhauser1978analysis}. Intuitively, submodularity of a set function describes a diminishing return property: adding an element to a smaller set gives a larger gain than adding it to a larger set. The formal definition is as follows.

\begin{definition}[\emph{Submodularity}] \label{De:Submodularity}
	A set function $f: 2^{\Omega} \rightarrow \mathbb{R}$ is called submodular if for all $A \subseteq B \subseteq \Omega$ and all elements $e \in \Omega,$ it holds that.
	\begin{equation} \label{Eq:SubmodularityDefinition}
			f(A \cup\{e\})-f(A) \geq f(B \cup\{e\})-f(B) .
	\end{equation}
\end{definition}
We find that submodularity does not hold for $J(S)$. 
Here we present a simple counterexample. Assume $\alpha_{1}=0.5, \alpha_{2}=2.5, \alpha_{3}=0.5$ and $\gamma_{s}=0.01, \gamma_{v}=0.05, \gamma_{u}=0.1$. Let $S_{1}=\{4,9,10\}$ and $S_{2}=\{2,3,4,9,10\}$, which implies $S_1 \subseteq S_2 $. Then we can compute directly that
$$
\begin{aligned} J\left(S_{1} \cup\{1\}\right)=-0.5982, J\left(S_{1}\right)=-0.5003; \\ J\left(S_{2} \cup\{1\}\right)=-0.7860, J\left(S_{2}\right)=-0.6910. \end{aligned}
$$
It is clear to see that
$$
\begin{aligned} &J\left(S_{1} \cup\{1\}\right)-J\left(S_{1}\right) =-0.098 \\ & \leq J\left(S_{2} \cup\{1\}\right)-J\left(S_{2}\right)=-0.095, \end{aligned}
$$
which violates \eqref{Eq:SubmodularityDefinition} in Definition \ref{De:Submodularity}, indicating that $J(S)$ is not submodular. Therefore, the greedy algorithm in previous works, e.g., \cite{summers2015submodularity}, cannot provide any guarantees when solving problem \eqref{Eq:H2ProblemFormulation}. 
Since our main focus is to identify the optimal formation of autonomous vehicles as shown in Problem~\ref{Pr:Main}, the brute force method is one straightforward approach to obtain the true optimal solution.



\section{NUMERICAL RESULTS OF OPTIMAL FORMATION}

\label{Section:Results}

In this section, we present the numerical results of optimal formation of multiple autonomous vehicles in mixed traffic flow based on formulation \eqref{Eq:H2ProblemFormulation}.

\subsection{Numerical Setup}
\begin{figure}[t]
	\centering
	\subfigure[]
	{ \label{Fig:OVMSpacing}
		\includegraphics[scale=0.4]{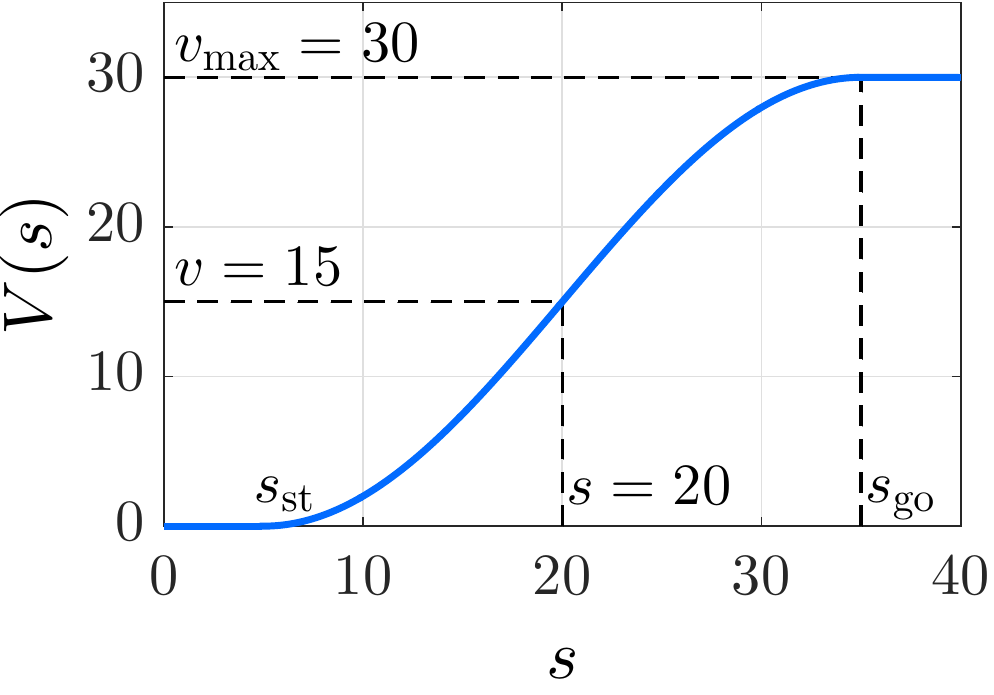}
	}
	\subfigure[]
	{\label{Fig:OVMVelocityDot}
		\includegraphics[scale=0.4]{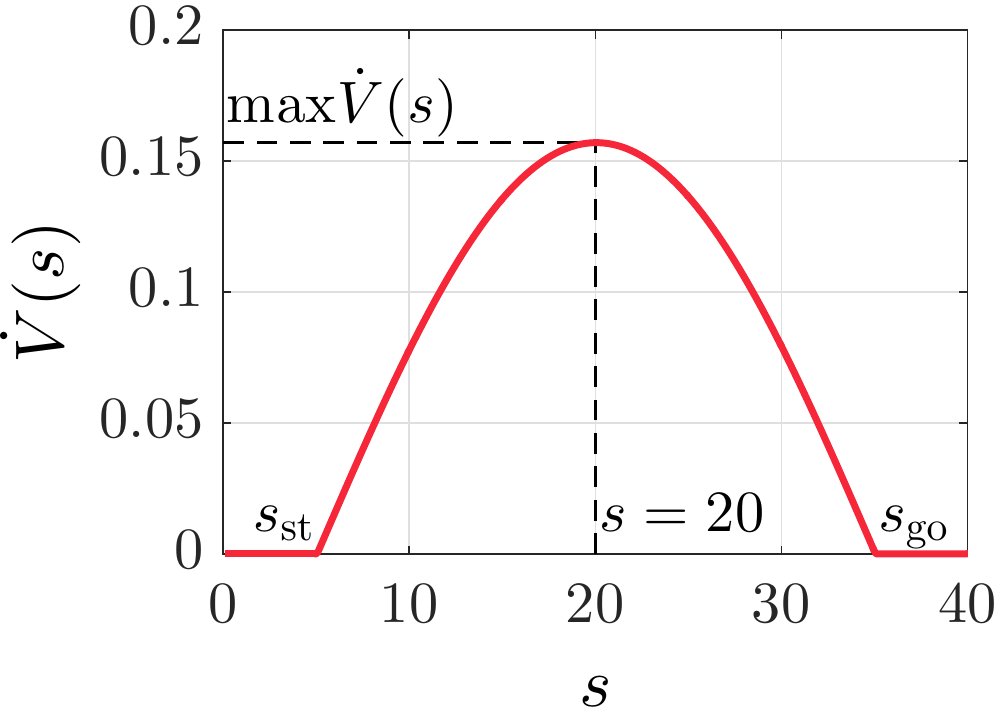}
	}
	\vspace{-3mm}
	\caption{Typical profile of the spacing-dependent desired velocity $V(s)$ and its derivative $ \dot{V}(s)$ when $\alpha = 0.6,\beta = 0.9, v_{\max} =30,s_{\text{st}}=5,s_{\text{go}}=35$.}
	\label{Fig:OVMPolicy}
\end{figure}

We consider an explicit car-following model, the optimal velocity model (OVM), for HDVs in our numerical studies. The specific model of \eqref{Eq:HDVNonlinearModel} can then be expressed as \citep{jin2016optimal}
\begin{equation}   \label{Eq:OVM}
F(\cdot)=\alpha \left(V(s_{i}(t))-v_{i}(t)\right)+\beta
\dot{s}_{i}(t),
\end{equation}
where $\alpha, \beta>0$ represent the driver's sensitivity coefficients and $V(\cdot)$ denotes the spacing-dependent desired velocity, typically given by
\begin{equation}
\label{Eq:OVM_DesiredVelocity}
V(s)=\begin{cases}
0, &s\le s_{\mathrm{st}};\\
f_v(s), &s_{\mathrm{st}}<s<s_{\mathrm{go}};\\
v_{\max}, &s\ge s_{\mathrm{go}},
\end{cases}
\end{equation}
with
\begin{equation}
\label{Eq:OVM_SpacingPolicy}
f_{v}(s)=\frac{v_{\max }}{2}\left(1-\cos (\pi
\frac{s-s_{\text{st}}}{s_{\text{go}}-s_{\text{st}}})\right).
\end{equation}
In OVM model, the values of the coefficients in the linearized HDV model \eqref{Eq:LinearHDVModel} can be calculated as
\begin{equation}
	\alpha_{1}=\alpha \dot{V}\left(s^{*}\right), \alpha_{2}=\alpha+\beta, \alpha_{3}=\beta,
\end{equation}
where $\dot{V}\left(s^{*}\right)$ denotes the derivative of $V(\cdot)$ at $s^*$. Fig.~\ref{Fig:OVMPolicy} illustrates the curves of $V(s)$ and $\dot{V}(s)$ under a typical parameter setup as that in \cite{jin2016optimal}.

\subsection{Two Predominant Optimal Formations}
The first numerical study aims to answer Problem \ref{Pr:Main}, i.e., identify the optimal formation of autonomous vehicles in mixed traffic flow. Formulation \eqref{Eq:H2ProblemFormulation} is considered and the brute force method is utilized for numerical computation. We fix $v_{\max} =30,s_{\text{st}}=5,s_{\text{go}}=35$ and let $\gamma_s=0.01,\gamma_v=0.05,\gamma_u=0.1$. Then we observe that the numerical solution of the optimal formation relies on the parameter setup in OVM model, i.e., the car-following behavior of HDVs. Three examples are listed in Table~\ref{Tb:OptimalFormation_Example} when $n=12,k=4$. Platoon formation, uniform distribution or even certain abnormal formations might be the optimal formation. A typical abnormal result is the same as that in Fig.~\ref{Fig:ExampleAbnormal}, which is essentially a transition pattern between platoon formation and uniform distribution.

\begin{table}[hb]
\begin{center}
\caption{Optimal Formation}\label{Tb:OptimalFormation_Example}
\vspace{-1mm}
\begin{tabular}{cccc}
$\alpha$ & $\beta$ & $s^*$  & numerical solution \\\hline
1.4 & 1.8 & 10  & platoon formation (Fig.~\ref{Fig:ExamplePlatoon})\\
0.6 & 0.9 & 20 & uniform distribution (Fig.~\ref{Fig:ExampleUniform})\\
0.9 & 1.3 & 16  & abnormal formation (Fig.~\ref{Fig:ExampleAbnormal})\\
 \hline
\end{tabular}
\end{center}
\end{table}

\begin{figure*}[t]
	\centering
	\subfigure[$n=12,k=2$]
	{ \label{Fig:FormationResult1}
		\includegraphics[scale=0.4]{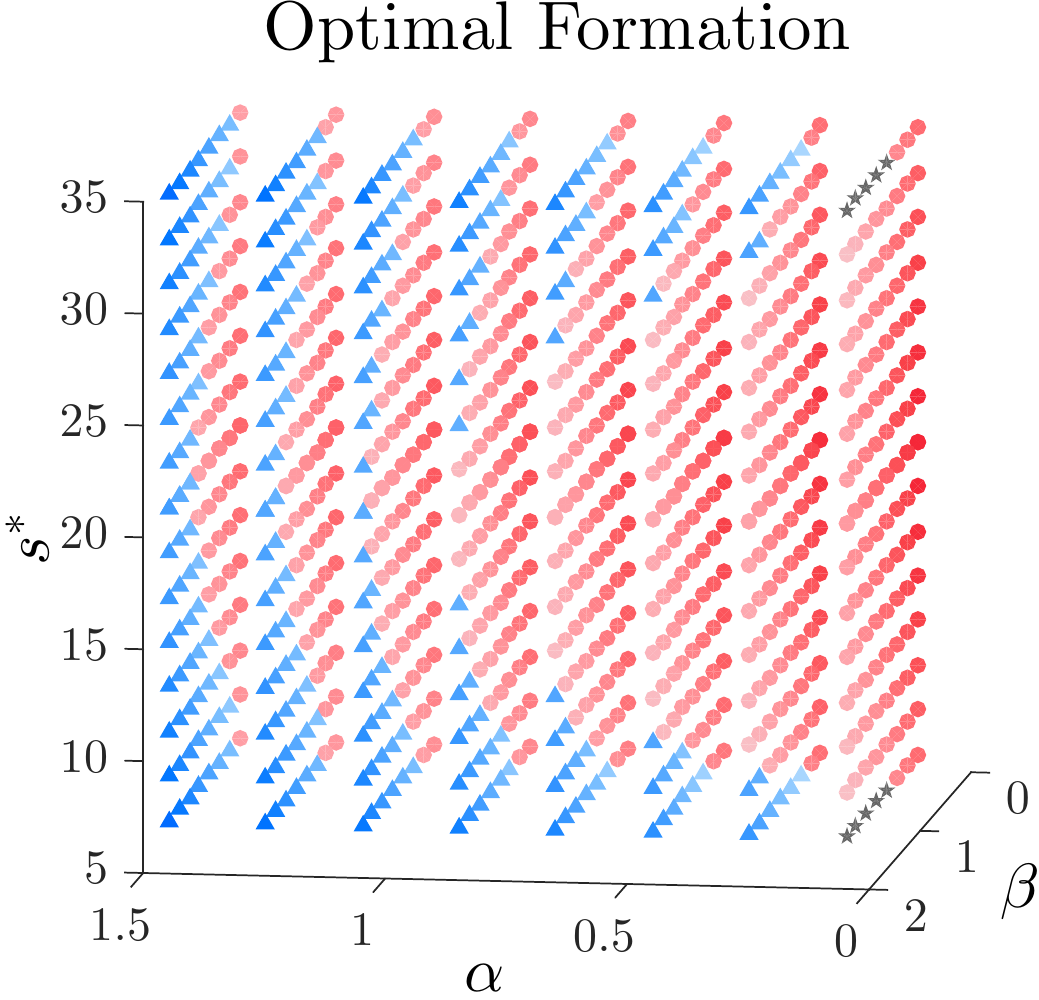}
		\includegraphics[scale=0.4]{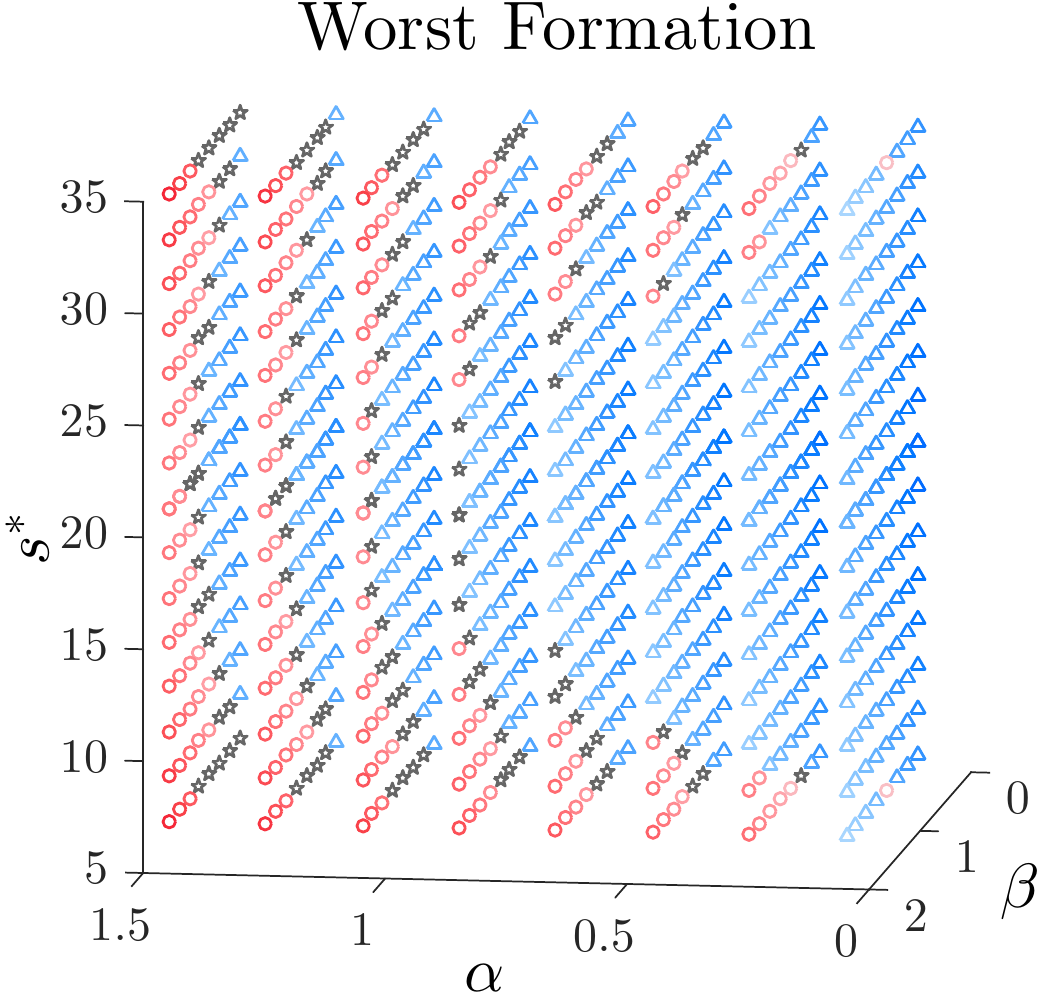}
	}
	\subfigure[$n=12,k=4$]
	{ \label{Fig:FormationResult2}
		\includegraphics[scale=0.4]{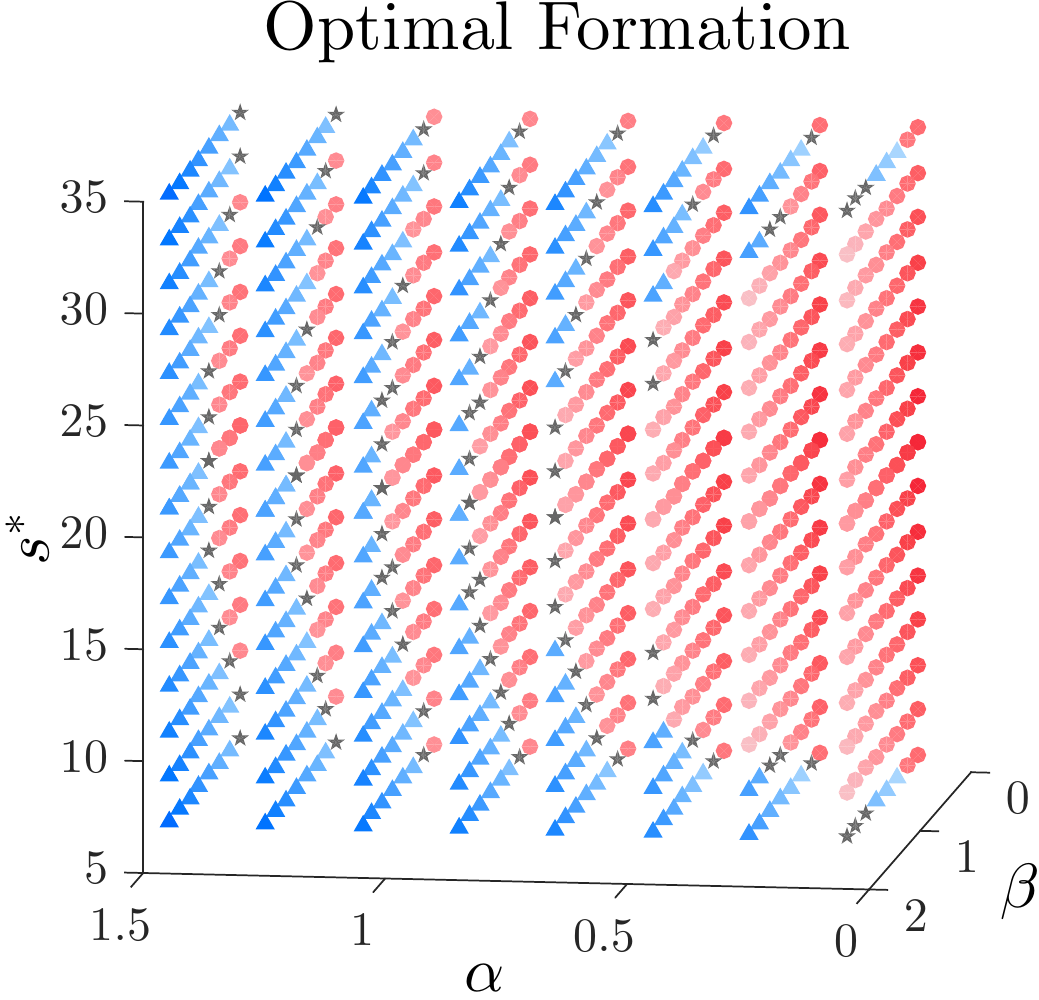}
		\includegraphics[scale=0.4]{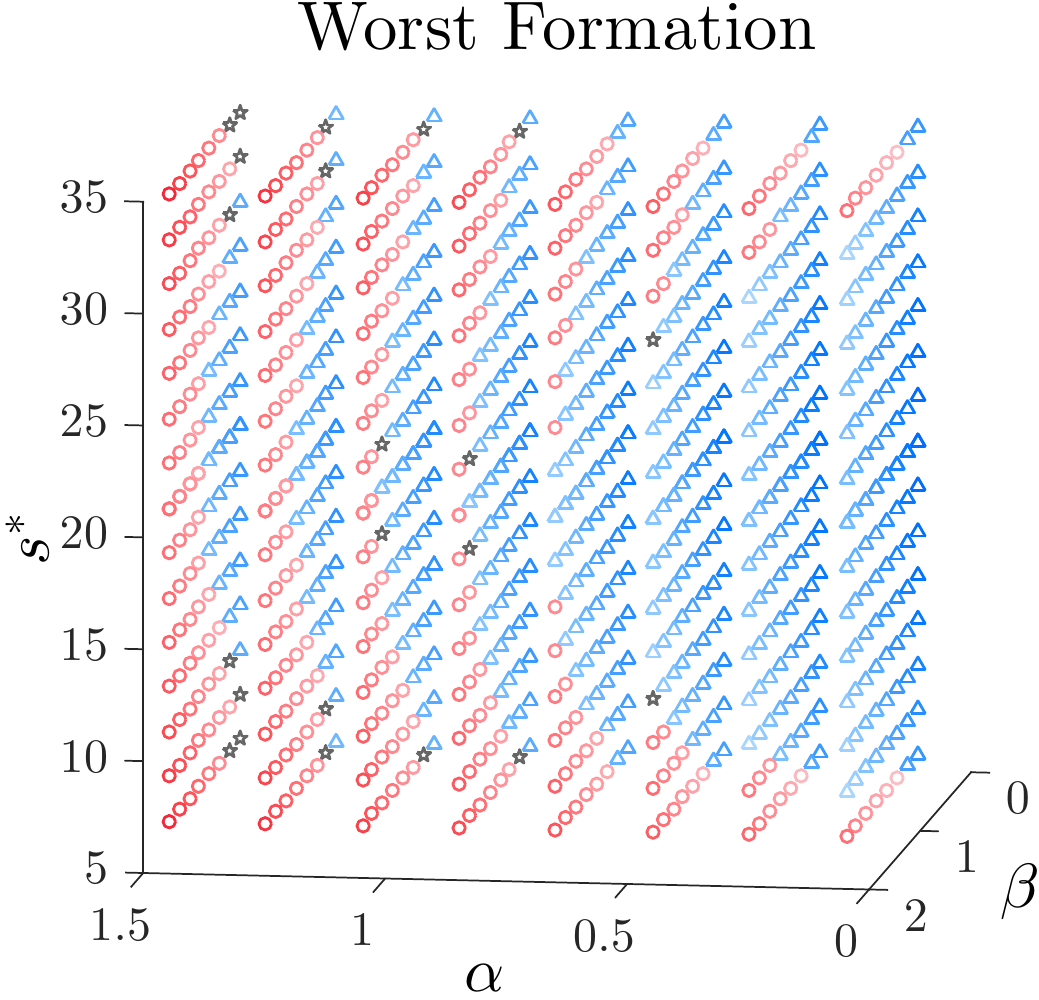}
	}
	\subfigure[$n=12,k=2$]
	{ \label{Fig:FormationResult3}
		\includegraphics[scale=0.4]{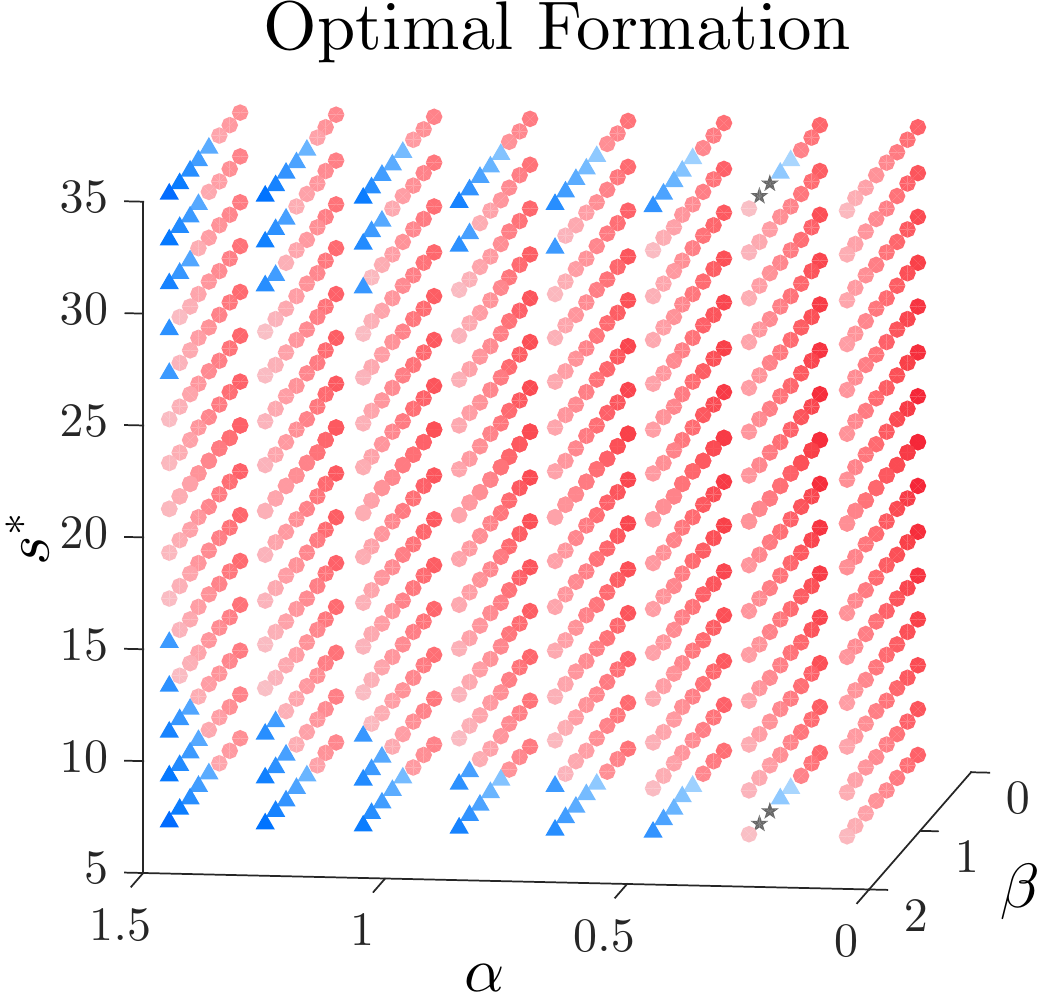}
		\includegraphics[scale=0.4]{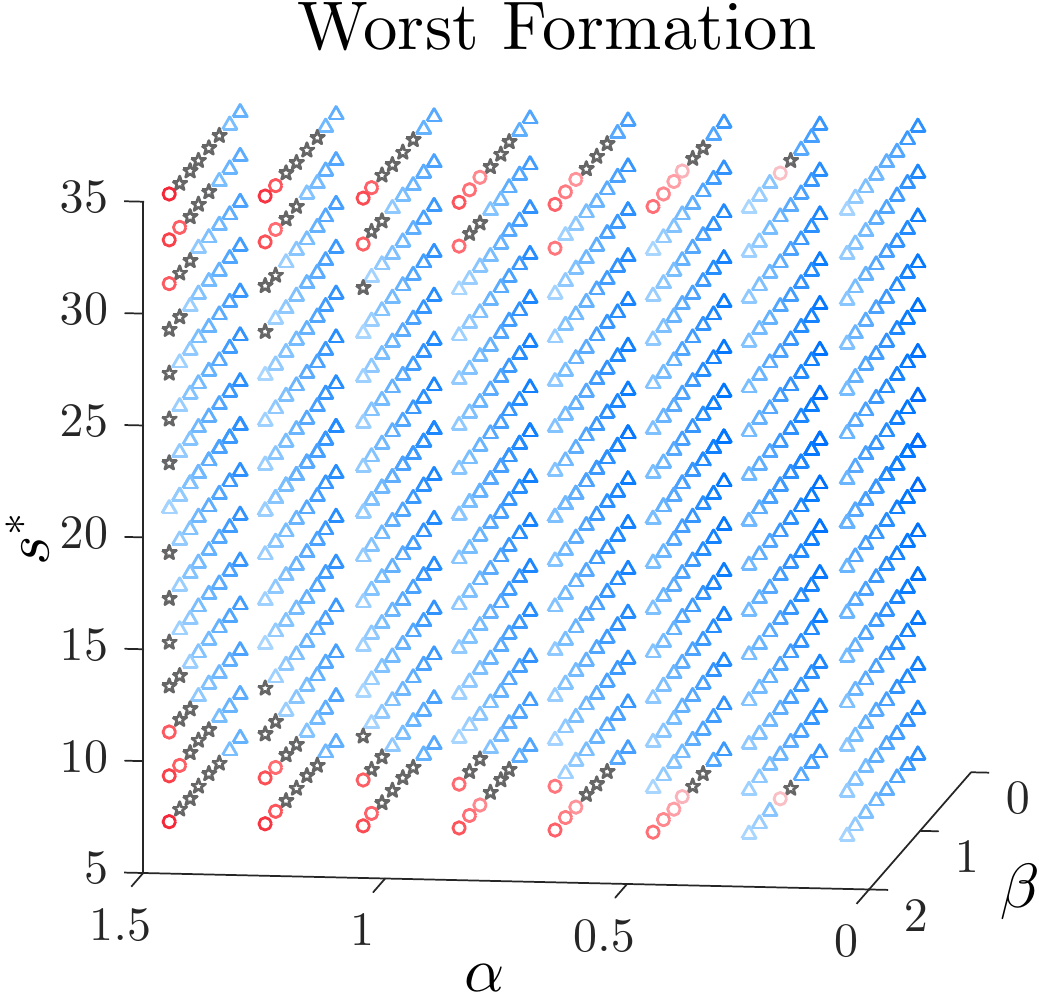}
	}
	\subfigure[$n=12,k=4$]
	{ \label{Fig:FormationResult4}
		\includegraphics[scale=0.4]{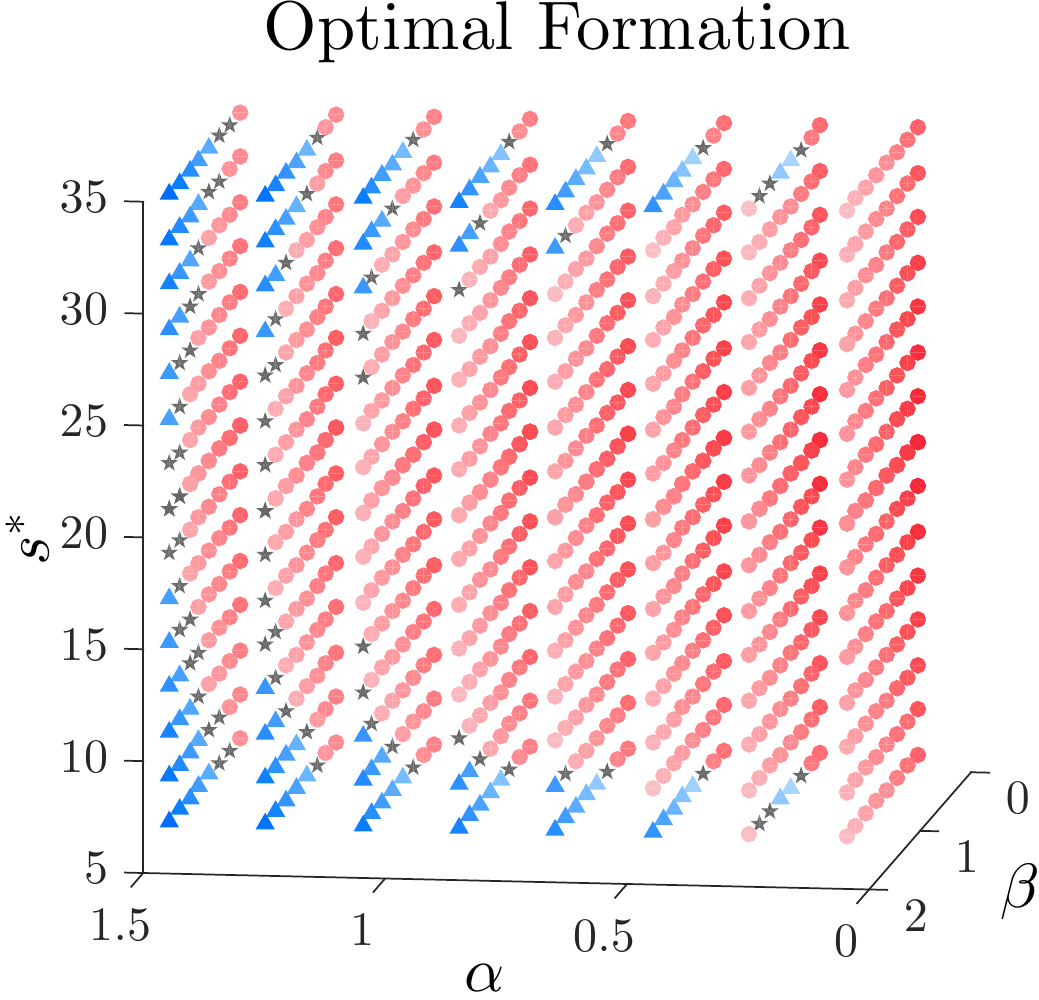}
		\includegraphics[scale=0.4]{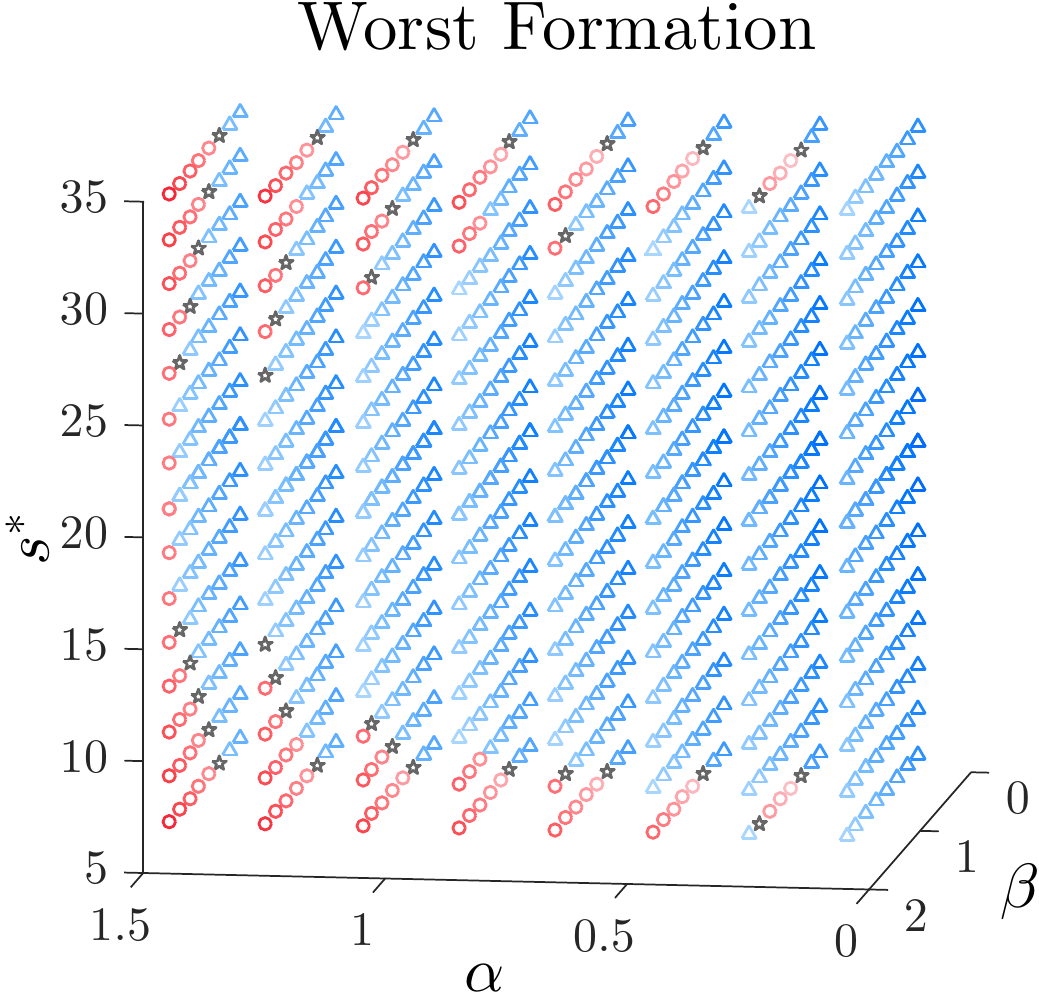}
	}
	\vspace{-3mm}
	\caption{Optimal and worst formation at various parameter setups. Red circles, blue triangles, and gray stars denote uniform distribution, platoon formation, and abnormal formations, respectively. In each panel, the left figure shows the optimal formation, where the darker the red, the larger the value of $\xi$; the darker the blue, the smaller the value of $\xi$. In contrast, the right figure shows the worst formation, where the darker the blue, the larger the value of $\xi$; the darker the red, the smaller the value of $\xi$. (a)(b) $\gamma_s=0.01,\gamma_v=0.05,\gamma_u=0.1$. (c)(d) $\gamma_s=0.03,\gamma_v=0.15,\gamma_u=0.1$.}
	\label{Fig:FormationResult}
\end{figure*}

We proceed to solve Problem \eqref{Eq:H2ProblemFormulation} in various parameter setups. The vehicle number is set to $n=12$, $k=2$ or $4$. We still fix $v_{\max} =30,s_{\text{st}}=5,s_{\text{go}}=35$, but discretize $\alpha,\beta,s^*$ within a common range: $\alpha \in [0.1,1.5],\beta\in [0.1, 1.5],s^*\in [5,35]$. Two different setups of the weight coefficients $\gamma_{s},\gamma_{v},\gamma_{u}$ in the performance output \eqref{Eq:Output} are also under consideration. Note that the brute force method can also offer the worst formation based on \eqref{Eq:H2ProblemFormulation} at the same time. The numerical results of optimal formation and worst formation are illustrated in Fig.~\ref{Fig:FormationResult}. As can be clearly observed, there exist two predominant patterns for optimal or worst formations: platoon formation and uniform distribution. This result holds regardless of the specific number of autonomous vehicles $k$ or the value of weight coefficients in \eqref{Eq:Output}. Between these two patterns there is an apparent boundary, along which there might exist some abnormal results. One interesting observation is that the optimal formation and the worst formation have an evident corresponding relationship. Precisely, when uniform distribution is optimal, platoon formation usually becomes the worst, and vice versa. This result indicates that the prevailing platoon formation might be the optimal formation, but could also be the worst choice, depending on the parameter of HDV car-following behaviors.

Then we make further investigations on the explicit relationship between the optimal formation and the HDV parameter setup. It is observed that the string stability performance of HDV car-following behaviors has a strong impact on the optimal formation of autonomous vehicles in mixed traffic flow. A series of vehicles is called string unstable if oscillations are amplified upstream the traffic flow. As shown in~\cite{orosz2010traffic}, the condition for strict string stability of OVM after linearization is
\begin{equation}
\alpha+2\beta \geq \dot{V} (s^*).
\end{equation}
Here we define a string stability index $\xi$ as
\begin{equation}
	\xi:=\alpha+2 \beta-\dot{V}\left(s^{*}\right).
\end{equation}
Note that a larger value of $\xi$ indicates better string stability
behavior. In our parameter setup, $\dot{V}\left(s^{*}\right)$ decreases as $\vert s^*-$
$20 \vert$ grows up, as shown in Fig.~\ref{Fig:OVMVelocityDot}. Therefore, a larger value of
$\alpha$, $\beta$ or $\vert s^*-20 \vert $ leads to a larger value of $\xi$, i.e., a better string stability performance of HDVs.

In Fig.~\ref{Fig:FormationResult}, we utilize the color darkness to indicate the value of $\xi$. Then the relationship between string stability of HDVs and the optimal formation of autonomous vehicles can be clearly observed. At a larger value of $\xi$, platoon formation could be the optimal choice. In contrast, when $\xi$ is small, indicating a poor string stability behavior of HDVs, uniform distribution achieves the best performance, while platoon formation becomes the worst. Note that in general cases we always assume that HDVs have a poor string stability behavior due to drivers' large reaction time and limited perception abilities \citep{sugiyama2008traffic}. This result indicates that platoon formation might have the most limited capability to improve traffic flow, compared to other possible formations in mixed traffic flow. One intuitive understanding is that when HDVs has poor string stability performance, distributing autonomous vehicles uniformly allows autonomous vehicles to maximize their capabilities in suppressing traffic instabilities and mitigating undesired perturbations. Instead, when all human drivers have better driving abilities, organizing all the autonomous vehicles into a platoon is a better choice.  

\subsection{Comparison Between Platoon Formation and Uniform Distribution}

\begin{figure}[tb]
	\centering
	\subfigure[]
	{ \label{Fig:SystemScaleComparison1}
		\includegraphics[scale=0.43]{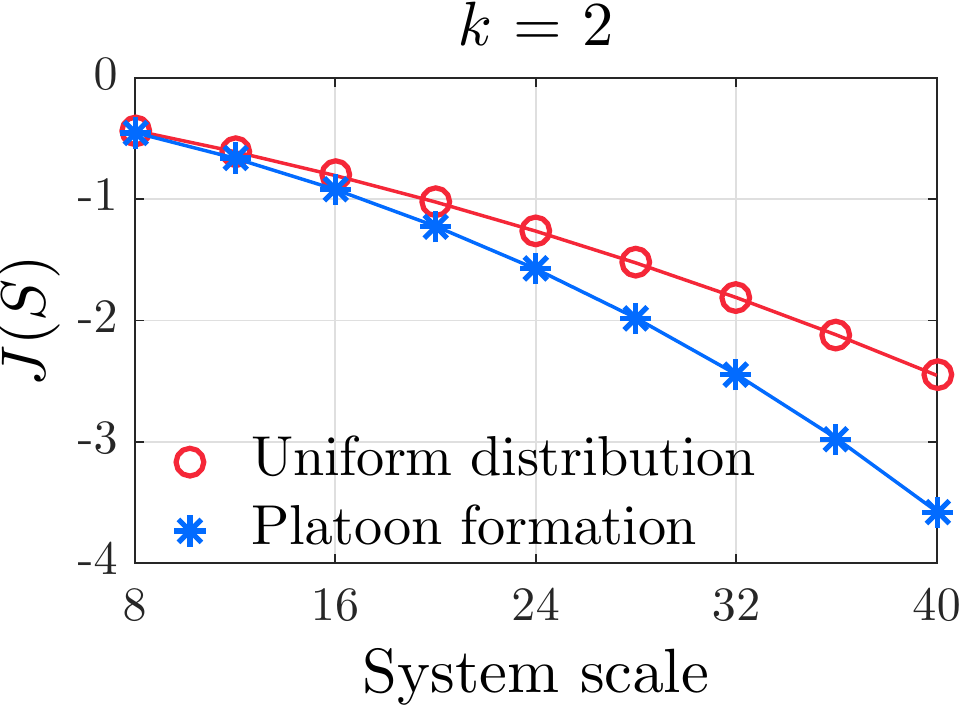}
		\includegraphics[scale=0.43]{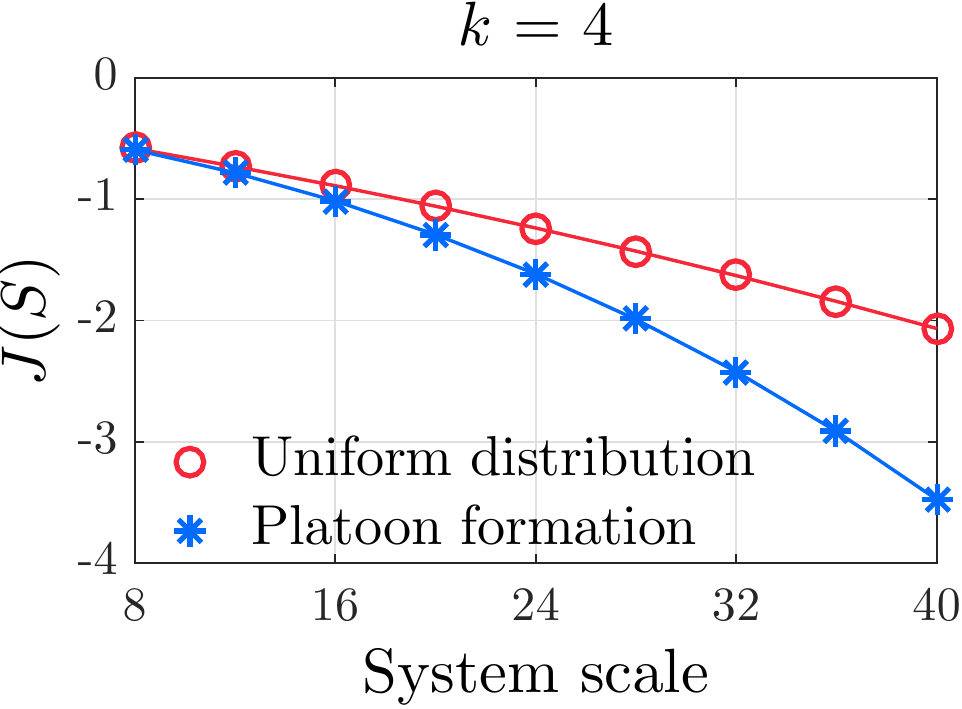}
	}
	\subfigure[]
	{\label{Fig:SystemScaleComparison2}
		\includegraphics[scale=0.43]{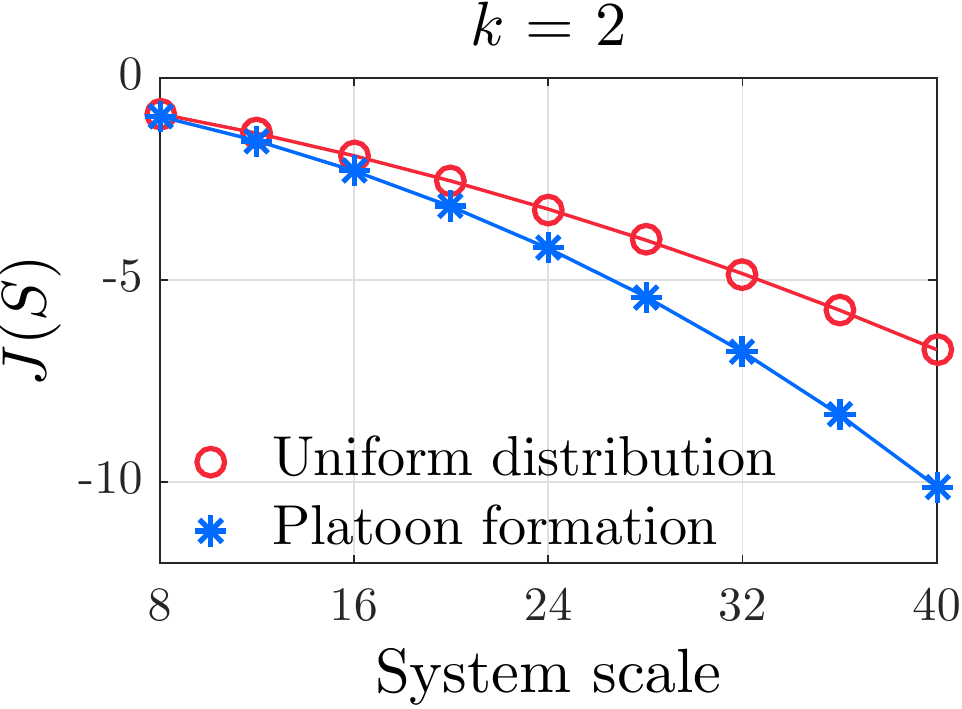}
		\includegraphics[scale=0.43]{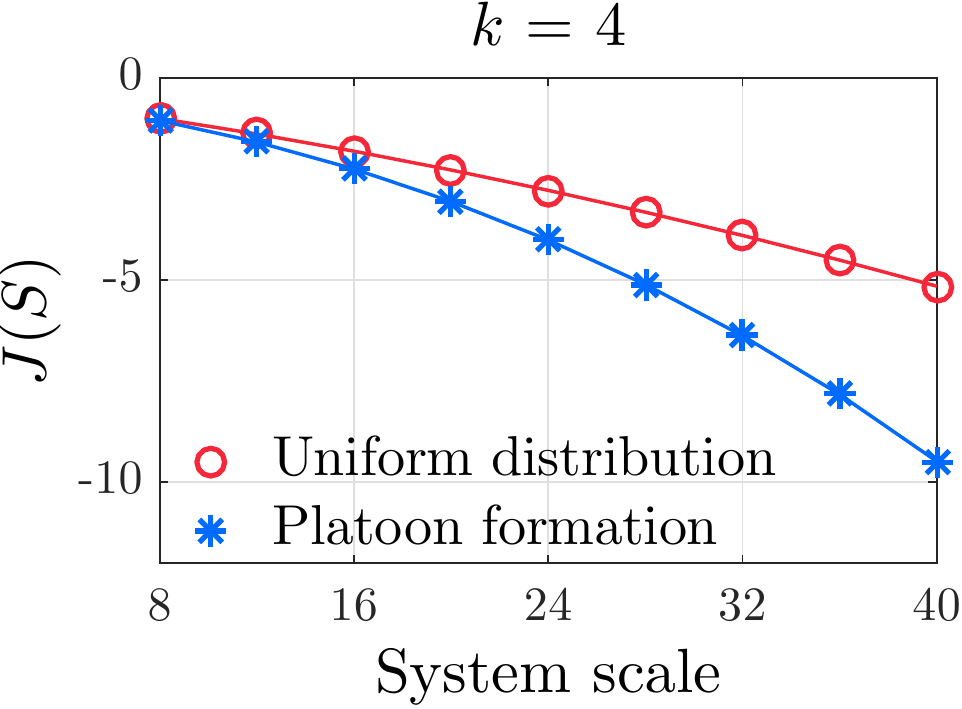}
	}
	\vspace{-3mm}
	\caption{Comparison between platoon formation and uniform distribution at different system scales. In OVM model,  $\alpha = 0.6,\beta = 0.9, s^*=20, v_{\max} =30,s_{\text{st}}=5,s_{\text{go}}=35$. (a) $\gamma_s=0.01,\gamma_v=0.05,\gamma_u=0.1$. (b) $\gamma_s=0.03,\gamma_v=0.15,\gamma_u=0.1$.}
	\label{Fig:SystemScaleComparison}
\end{figure}

Finally, we carry out another numerical study to make further comparisons between the two predominant formations at different system scales $n \in [8,40]$. In Section 4.2, we consider different OVM parameters, but focus on the single case where $n=12$. Here we vary the system scale, but fix the OVM model to a typical setup for human's driving behavior as that in \cite{jin2016optimal}. The comparison results for $J(S)$ of these two formations are demonstrated in Fig.~\ref{Fig:SystemScaleComparison} ($k=2$ or $4$). Recall that a larger value of $J(S)$ denotes a better performance when $\vert S \vert$ is fixed. We observe that in this parameter setup, uniform distribution is optimal, while platoon formation is the worst. Moreover, as shown in Fig.~\ref{Fig:SystemScaleComparison}, the performance gap between the two formations is rapidly enlarged as the system scale grows up. This results indicates that at a large system scale, i.e., a low penetration rate of autonomous vehicles, there could exist a huge performance difference between platoon formation and other possible formations, e.g., uniform distribution. In the near future when we only have a few autonomous vehicles on the road, platooning might not be a suitable choice for improving the global traffic performance.

\section{Conclusion}

\label{Section:Conclusion}

In this paper, we have formulated a set function optimization problem to investigate the optimal formation for autonomous vehicles in mixed traffic flow. Taking into account the $\mathcal{H}_2$ performance metric, we reveal 
that there exist two predominant optimal formations for autonomous vehicles: uniform distribution and platoon formation. In particular, our results indicate that when HDVs have a poor string stability behavior, the prevailing vehicle platooning is not a suitable choice, which might even have the least potential in mitigating traffic perturbations. 
Our results suggest more opportunities for the formation of autonomous vehicles in mixed traffic flow, beyond the prevailing platoon formation. Reconsidering these possibilities can take full advantage of the benefits brought by autonomous vehicles in mixed traffic systems.


\bibliography{ifacconf}             








\end{document}